\documentclass[11pt,leqno]{amsart}

\usepackage[letterpaper,hcentering,textwidth=6in,textheight=8.55in,headsep=15pt]{geometry}
\usepackage{amsmath,amssymb}
\usepackage{mathptmx}
\DeclareSymbolFontAlphabet{\ptmcal}{symbols}
\DeclareSymbolFont{cmcal}{OMS}{cmsy}{m}{n}
\DeclareSymbolFontAlphabet{\mathcal}{cmcal}
\usepackage[hidelinks]{hyperref}

\allowdisplaybreaks

\newtheorem{theorem}{Theorem}[section]
\newtheorem{lemma}[theorem]{Lemma}
\newtheorem{proposition}[theorem]{Proposition}
\numberwithin{equation}{section}

\newcommand{\papertitle}{
Near-Diagonal Lower Bounds Do Imply \\ Sub-Gaussian Upper Bounds for Heat Kernels:\\ Where Analysis Meets Probability
}
\title[Heat Kernel Upper Bounds]{\papertitle}
\author{Bo Li}
\address{ Institute of Mathematics, Jiaxing University, Jiaxing 314001, P. R. China.}
\email{bli@zjxu.edu.cn}
\author{Yuhua Sun}
\address{School of Mathematical Sciences and LPMC, Nankai University,
	Tianjin 300071, P. R. China}
\email{sunyuhua@nankai.edu.cn}

\hypersetup{pdftitle={\papertitle},pdfsubject={Grigor'yan--Hu--Lau Conjecture 3.6}}

\makeatletter
\@ifundefined{subjclassname@2020}{\@namedef{subjclassname@2020}{\textup{2020} Mathematics Subject Classification}}{}
\def\@settitle{\begin{center}\parbox{0.84\textwidth}{\centering\normalfont\Large\bfseries\@title}\end{center}}
\makeatother

\subjclass[2020]{Primary 35K08, 31C25; Secondary 60J35.}
\keywords{Heat kernel, Dirichlet form, near-diagonal lower bound, exit time, sub-Gaussian estimate}
\begin{document}

\begin{abstract}
We prove that a near-diagonal lower bound for the heat kernel of a local
regular Dirichlet form on a metric measure space with a lower
$\alpha$-regular measure implies a full sub-Gaussian upper bound.
This settles a conjecture posed by Grigor'yan, Hu, and Lau in
\emph{Comm. Pure Appl. Math.} (2008).
We provide two proofs. One is analytic, and the other is probabilistic.
In the analytic proof, we use locality to compare resolvent potentials
across separating annuli and iterate the comparison to obtain
off-diagonal resolvent decay.
In the probabilistic proof, we sample continuous paths on deterministic
time grids and use the Markov property to derive an exponential bound
on the probability of reaching a prescribed radius within a short time.
\end{abstract}

\maketitle
\markboth{HEAT KERNEL UPPER BOUNDS}{HEAT KERNEL UPPER BOUNDS}

\section{Introduction and the main result}

Upper and lower heat kernel bounds describe the spatial and temporal behavior of the associated Markov processes.
Classical heat kernel estimates have been well established in the following settings:
\begin{itemize}
  \item complete Riemannian manifolds, including the Li--Yau estimates and Gaussian upper and lower bounds \cite{cy81,davies87b,g97,ly86,sc92}.

  \item local and non-local Dirichlet forms on metric measure spaces, including Gaussian and sub-Gaussian estimates, on/off-diagonal bounds, and gradient estimates \cite{ckw21,cjks20,gh08,ghh17,kkk22,sturm96}.

  \item Riemannian curvature-dimension (RCD) spaces, including Gaussian upper/lower bounds for the heat kernel and corresponding gradient estimates  \cite{jiang15,jlz16}.

  \item Davies--Gaffney estimates and finite propagation methods for abstract heat semigroups and self-adjoint operators \cite{cgt82,cs08}.

  \item Euclidean domains with Dirichlet boundary conditions, including bounded and exterior domains \cite{chenkimSong10,chenkimSong16,qizhang03}.

  \item special geometric, discrete, and fractal structures,
  including  the unit sphere \cite{nss19}, hyperbolic spaces and Kleinian groups \cite{daviesmand88},
  Markov chains \cite{cg97}, fractals \cite{bp88}, and strip complexes \cite{bendikov11}.
\end{itemize}
For further background on heat kernel estimates and spectral theory, 
we refer to the monographs of Davies \cite{davies89} and Grigor'yan \cite{grigor09}.

A common approach to two-sided estimates, going back to Aronson \cite{aronson67}, is to establish the upper bound first.
The upper bound controls the heat kernel mass away from the starting point.
A lower bound for the mass inside a ball then gives an on-diagonal lower estimate
through the semigroup identity and the Cauchy--Schwarz inequality.
Under suitable hypotheses, regularity estimates or a chaining argument extend the lower bound away from the diagonal.
In 2008, Grigor'yan, Hu, and Lau \cite{ghl08} \textit{reversed this order}
by deriving an on-diagonal upper bound from a global near-diagonal lower bound and a lower volume bound.
They described this as the first result of its kind and conjectured that
the same hypotheses imply a full sub-Gaussian upper bound for local regular Dirichlet forms \cite[Conjecture~3.6]{ghl08}.
In the same paper, Grigor'yan, Hu, and Lau showed that
a stronger local lower estimate implies a full sub-Gaussian upper bound \cite[Theorem~4.2]{ghl08}.
They also proved that a near-diagonal lower estimate, together with a time-independent upper estimate,
yields the full upper bound \cite[Theorem~4.6]{ghl08}.
Subsequent work identified additional conditions that give the off-diagonal bound.
On doubling spaces with uniform volume growth,
Boutayeb \cite[Theorem~1.1]{boutayeb09} showed that a global
sub-Gaussian lower bound implies an on-diagonal upper bound.
For conservative local regular Dirichlet forms on doubling spaces,
he obtained a full sub-Gaussian upper bound from a global sub-Gaussian
lower bound under a two-sided mean exit time estimate
\cite[Theorem~1.2]{boutayeb09}.
For non-local regular Dirichlet forms without a killing part on
Ahlfors regular metric measure spaces,
Grigor'yan, Hu and Hu \cite[Theorem~2.9]{ghh17}
showed that the near-diagonal lower estimate,
two-sided bounds on the jump kernel,
and a survival estimate together are equivalent to two-sided
stable-like heat kernel bounds.
Despite these developments, to the best of our knowledge,
the conjecture of Grigor'yan, Hu, and Lau in \cite{ghl08} has remained open under its original assumptions.
In this paper, we completely resolve their Conjecture 3.6 by proving that a near-diagonal lower bound,
together with lower volume regularity,
implies the full sub-Gaussian upper bound for local regular Dirichlet forms, \textit{without} any additional assumptions.

Let $(M,d)$ be a locally compact separable metric space with a Radon measure $\mu$ of full support.
Write $B(x,r)=\{y\in M:d(x,y)<r\}$ for the open ball,
and $V(x,r)=\mu(B(x,r))$ for its volume.
Let $(\ptmcal{E},\ptmcal{F})$ be a local regular symmetric Dirichlet form on $L^2(M,\mu)$ with heat semigroup $(P_t)_{t>0}$ and heat kernel $p_t(x,y)$. For each fixed time, kernel estimates are understood to hold for $\mu\otimes\mu$-almost every pair.
For $0< r_0\le\infty$, $\alpha>0$, and $\beta>1$, we say that $\mu$ is lower $\alpha$-regular,
and $p_t(x,y)$  satisfies a near-diagonal lower estimate/an on-diagonal upper estimate/an off-diagonal upper estimate, respectively,
if there are $c_1,c_2,c_3,c_4,C_0>0$ and $0<\delta\le1/2$ such that
\begin{align}
  V(x,r)&\ge c_1 r^\alpha,
  &&0<r<r_0,
  \tag{\ensuremath{V_{\ge}}}\label{eq:vol}\\
  p_t(x,y)&\ge c_2 t^{-\alpha/\beta},
  &&0<t<\delta r_0^\beta,\ d(x,y)<\delta t^{1/\beta}
  \tag{\ensuremath{NLE}}\label{eq:nle}\\
  p_t(x,y)&\le c_3 t^{-\alpha/\beta},
  &&0<t<\delta r_0^\beta,
  \tag{\ensuremath{DUE}}\label{eq:due}\\
   p_t(x,y)&\le C_0 t^{-\alpha/\beta}
  \exp\left[-c_4\left(\frac{d(x,y)}{t^{1/\beta}}\right)^{\beta/(\beta-1)}\right],
  &&0<t< r_0^\beta. \tag{\ensuremath{UE}}\label{eq:ue}
\end{align}
The volume bound is required for $\mu$-almost every center.

\begin{theorem}\label{thm:main}
Let \((M, d, \mu)\) be a locally compact separable metric measure space endowed with a local regular symmetric Dirichlet form 
 \((\ptmcal{E}, \ptmcal{F})\) in \(L^2(M, \mu)\).
If \(\mu\) admits \eqref{eq:vol}, then
\[\eqref{eq:nle}
\implies \eqref{eq:ue}.
\]
\end{theorem}

We remark that it suffices to prove the bound in \eqref{eq:ue} for
$0<t<\delta r_0^\beta$.
Indeed, suppose that the bound holds for $0<t<T$,
and fix $0<t<2T$.
In the semigroup formula for $p_t(x,y)$,
each integration point $z$ satisfies
$d(x,z)\geq d(x,y)/2$ or $d(z,y)\geq d(x,y)/2$.
On each of these two regions,
we apply the upper bound at time $t/2$ to the corresponding
kernel factor and integrate the other factor using symmetry
and the sub-Markov property
to obtain
\[
p_t(x,y)
\leq 2C_0(t/2)^{-\alpha/\beta}
\exp\left[
-c_4\left(\frac{d(x,y)/2}{(t/2)^{1/\beta}}\right)^{\beta/(\beta-1)}
\right].
\]
Thus the same type of bound holds for $0<t<2T$.
If $r_0<\infty$,
finitely many iterations starting from $T=\delta r_0^\beta$
extend the bound to $0<t<r_0^\beta$,
with adjusted constants.
Moreover,  \eqref{eq:nle} and the sub-Markov property imply that
$V(x,r)\leq C r^\alpha$ for $0<r<\delta^{1+1/\beta}r_0$.
Together with \eqref{eq:vol},
this shows that $\mu$ is Ahlfors $\alpha$-regular on this range,
and on all scales if $r_0=\infty$.
The main result can be extend to globally doubling measures in
volume-normalized form, with $t^{-\alpha/\beta}$ replaced by
$V(x,t^{1/\beta})^{-1}$ in both \eqref{eq:nle} and \eqref{eq:ue}.
Volume doubling gives the uniform lower bound on the ball mass
needed in the annulus argument,
and a weighted Nash estimate yields the corresponding on-diagonal upper bound.

To show Theorem~\ref{thm:main}, we provide two proofs.
One is analytic, and the other is probabilistic.
Both use \eqref{eq:vol} and \eqref{eq:nle} to obtain a uniform lower bound
for the heat kernel mass in an annulus when the starting point
lies in a narrower concentric annulus.
In the analytic proof in Section~\ref{sec:upper-analytic},
we integrate this lower bound in time
and use a comparison principle from locality to transfer bounds
from the annuli to an inner ball.
Summing over the disjoint annuli and iterating the resulting estimate
give exponential resolvent decay,
which we transfer to the heat semigroup.
In the probabilistic proof in Section~\ref{sec:upper-prob},
path continuity forces a diffusion reaching the outer radius
before the end of its lifetime to cross each annulus.
Conditional on a visit sampled at a grid time,
the lower bound gives a uniform probability of finding the process
in the corresponding annulus at a common later time.
Summing these probabilities bounds the probability of rapid passage,
and iteration on deterministic time grids yields an exponential passage-time bound.
Both proofs combine the resulting semigroup estimates
with the on-diagonal upper estimate of Grigor'yan, Hu, and Lau
to obtain the full sub-Gaussian upper bound for the heat kernel.

The letters $C$ and $c$ denote positive constants whose values may change between estimates.

\section{Sub-Gaussian Upper Bounds via Analytic Methods}
\label{sec:upper-analytic}

The analytic proof is based on a resolvent estimate.  The key step
is a comparison principle arising from locality, which allows us
to propagate an upper bound across a separating annulus.
We then convert the resolvent estimate into an off-diagonal
estimate for the heat semigroup and combine it with the
on-diagonal upper estimate.

\subsection{A comparison principle from locality}

For an open set $O\subset M$, let $\ptmcal{F}_O$ denote the
domain of the part Dirichlet form on $O$, with its elements extended by zero to $M$.
For $\lambda>0$, write
\[
  \ptmcal{E}_\lambda(u,v)
  =
  \ptmcal{E}(u,v)+\lambda\langle u,v\rangle.
\]

\begin{lemma}
\label{lem:analytic-comparison}
Let $O_1,O_2\subset M$ be disjoint open sets, and let $U$
be an open neighborhood of $M\setminus(O_1\cup O_2)$.
Suppose that $g\in\ptmcal{F}\cap L^\infty(M)$ is real-valued and
$g\leq 0$ almost everywhere on $U$.
If $\ptmcal{E}_\lambda(g,v)\leq 0$ for every nonnegative $v\in\ptmcal{F}_{O_1}$, then
$g\leq 0 $ almost everywhere on $O_1$.
\end{lemma}

\begin{proof}
Put $G=O_1\cup O_2$.
Since $g\leq0$ almost everywhere on the open neighborhood
$U$ of $M\setminus G$, the quasi-continuous representative
of $g^+$ vanishes quasi-everywhere on $M\setminus G$.
Consequently, $g^+\in\ptmcal{F}_G$.
By regularity of the part form on $G$ and locality, the
decomposition over the disjoint open sets $O_1$ and $O_2$
gives $  g^+\mathbf{1}_{O_1}\in\ptmcal{F}_{O_1}$.
Taking $v=g^+\mathbf{1}_{O_1}$
as a test function,
we have $v\geq0$ and $\ptmcal{E}_\lambda(g,v)\leq0$.
By locality and the decomposition over $O_1$ and $O_2$,
we have $\ptmcal{E}_\lambda(g^+\mathbf{1}_{O_2},v)=0$.
Since $g^-v=0$ almost everywhere,
the contraction property gives
$\ptmcal{E}_\lambda(g^-,v)\leq0$.
It follows that
\[
0
\geq \ptmcal{E}_\lambda(g,v)
= \ptmcal{E}_\lambda(v,v)-\ptmcal{E}_\lambda(g^-,v)
\geq \ptmcal{E}_\lambda(v,v)
\geq \lambda\|v\|_2^2,
\]
which implies $v=0$ and hence $g\leq 0$ almost everywhere on $O_1$.
\end{proof}

\subsection{Resolvent and semigroup estimates}

Let $L$ be the nonnegative self-adjoint operator associated
with $(\ptmcal{E},\ptmcal{F})$, so that $P_t=e^{-tL}$.
For $\lambda>0$, define the resolvent $  R_\lambda=(\lambda+L)^{-1}$.
By the spectral theorem,
\begin{equation}\label{eq:analytic-resolvent-representation}
  R_\lambda f
  =
  \int_0^\infty e^{-\lambda s}P_s f\,\mathrm{d}s,
  \qquad f\in L^2(M,\mu).
\end{equation}
The identity holds in the strong operator sense on
$L^2(M,\mu)$ and, for nonnegative measurable functions,
by monotone approximation.
The resolvent equation is
\[
  \ptmcal{E}_\lambda(R_\lambda f,\varphi)
  =
  \langle f,\varphi\rangle,
  \qquad f\in L^2(M,\mu),\quad \varphi\in\ptmcal{F}.
\]

\begin{lemma}\label{lem:resolvent-decay}
There are $C,c,\lambda_0>0$ such that, for every Borel set
$F\subset M$ and almost every $x\in M$,
\[
  R_\lambda\mathbf{1}_F(x)
  \leq \frac{C}{\lambda}
  \exp\left[-c\,d(x,F)\lambda^{1/\beta}\right],\qquad \lambda\geq\frac{\lambda_0}{\delta r_0^\beta}.
\]
\end{lemma}

\begin{proof}
Without loss of generality, we assume that $d(x,F)>0$.
Integrating \eqref{eq:nle} and using \eqref{eq:vol},
we obtain for each fixed $0<s<\delta r_0^\beta$ and almost every $z\in M$
\begin{align*}
0<c_1c_2\delta^\alpha\le c_2s^{-\alpha/\beta}V(z,\delta s^{1/\beta})
&\leq \int_{B(z,\delta s^{1/\beta})}p_s(z,y)\,\mu(\mathrm{d}y)
\leq 1,
\end{align*}
where the upper bound follows from the sub-Markov property.
First choose a sufficiently large integer $N$ and then a constant
$\lambda_0>0$, so that
\[
c_1c_2\delta^\alpha N\geq 4,
\qquad
\lambda_0>\left(\frac{2(4N+2)}{\delta}\right)^\beta.
\]
Fix $\lambda>0$ with $\lambda\geq\lambda_0/(\delta r_0^\beta)$,
where the lower threshold is understood to be zero if $r_0=\infty$.
Every ball $B(x,(4N+2)\lambda^{-1/\beta})$ has finite measure.
Indeed, at the time
\[
s=\frac{1}{\lambda}\left(\frac{2(4N+2)}{\delta}\right)^\beta
<\delta r_0^\beta
\]
the triangle inequality gives
\[
B(x,(4N+2)\lambda^{-1/\beta})
\subset B(z,\delta s^{1/\beta}),
\]
for some $z\in B(x,(4N+2)\lambda^{-1/\beta})$ for which the volume
upper bound at time $s$ holds, and the latter ball has finite measure.

For $x\in M$ and $1\leq j\leq N$,
define the open annuli centered at $x$ by
\begin{align*}
  A_j &= \{y\in M: (4j-2)\lambda^{-1/\beta}<d(x,y)<(4j+2)\lambda^{-1/\beta}\},\\
  B_j &= \{y\in M: (4j-1)\lambda^{-1/\beta}<d(x,y)<(4j+1)\lambda^{-1/\beta}\}.
\end{align*}
The annuli $A_j$ are pairwise disjoint.
Each $A_j$ has finite measure.
If $z\in B_j$ and $0<s<1/\lambda$,
the triangle inequality gives $B(z,\delta s^{1/\beta})\subset A_j$.
Hence $P_s\mathbf{1}_{A_j}(z)\geq c_1c_2\delta^\alpha$
for almost every $z\in B_j$ at each such $s$.
By \eqref{eq:analytic-resolvent-representation}, we obtain that for almost every $z\in B_j$,
\[
R_\lambda\mathbf{1}_{A_j}(z)
\geq c_1c_2\delta^\alpha\int_0^{1/\lambda}e^{-\lambda s}\,\mathrm{d}s
=\frac{c_1c_2\delta^\alpha(1-e^{-1})}{\lambda}
\geq \frac{2}{N\lambda}.
\]

Suppose first that $F$ is nonempty and has finite measure.
We prove by induction on $k\geq0$ that
\[
R_\lambda\mathbf{1}_F(x)\leq\frac{1}{2^k\lambda}
\quad\text{for almost every }x\text{ with }\
d(x,F)>k(4N+2)\lambda^{-1/\beta}.
\]
The case $k=0$ follows from
$R_\lambda\mathbf{1}_F\leq R_\lambda\mathbf{1}_M\leq\lambda^{-1}$.
Suppose the assertion holds for $k-1$,
and fix $x$ with $d(x,F)>k(4N+2)\lambda^{-1/\beta}$.
For almost every $z\in B(x,(4N+2)\lambda^{-1/\beta})$,
the induction hypothesis gives
$R_\lambda\mathbf{1}_F(z)\leq (2^{k-1}\lambda)^{-1}$,
and $\mathbf{1}_F(z)=0$.
For each $1\leq j\leq N$,
use the annuli centered at $x$ and put
\[
g=R_\lambda\mathbf{1}_F-\frac{N}{2^k}R_\lambda\mathbf{1}_{A_j}.
\]
Both indicators belong to $L^2(M,\mu)$,
and the resolvent bounds give $g\in\ptmcal{F}\cap L^\infty(M)$.
The lower bound on $R_\lambda\mathbf{1}_{A_j}$ shows that
\[
g(z)
=R_\lambda\mathbf1_F(z)
-\frac{N}{2^k}R_\lambda\mathbf1_{A_j}(z)\leq \frac{1}{2^{k-1}\lambda}
-\frac{N}{2^k}\frac{2}{N\lambda}
=0\quad \text{for almost everywhere } z\in B_j.
\]
Apply Lemma~\ref{lem:analytic-comparison} with
\[
O_1=B(x,4j\lambda^{-1/\beta}),\quad
O_2=\{y\in M:d(x,y)>4j\lambda^{-1/\beta}\},
\quad U=B_j.
\]
Indeed, $M\setminus(O_1\cup O_2)\subset B_j$,
and for every nonnegative $v\in\ptmcal{F}_{O_1}$,
the resolvent equation gives
\[
\ptmcal{E}_\lambda(g,v)
=\langle\mathbf{1}_F,v\rangle
-\frac{N}{2^k}\int_{A_j}v\,\mathrm{d}\mu
\leq 0.
\]
We conclude that
$R_\lambda\mathbf{1}_F(z)\leq (N/2^k)R_\lambda\mathbf{1}_{A_j}(z)$
for almost every $z\in B(x,\lambda^{-1/\beta})\subset O_1$.
We sum over $1\leq j\leq N$ and use the disjointness of the annuli $A_j$ to obtain
\[
N R_\lambda\mathbf{1}_F(z)
\leq \frac{N}{2^k}\sum_{j=1}^N R_\lambda\mathbf{1}_{A_j}(z)
\leq \frac{N}{2^k}R_\lambda\mathbf{1}_M(z)
\leq \frac{N}{2^k\lambda}\quad \text{for almost every }z\in B(x,\lambda^{-1/\beta}).
\]
As $x$ varies,
these balls cover $\{z\in M:d(z,F)>k(4N+2)\lambda^{-1/\beta}\}$.
Taking a countable subcover completes the induction.

The exceptional sets may be taken simultaneously for all $k\geq 0$.
Choose the largest $n\geq 0$ with
$n(4N+2)\lambda^{-1/\beta}<d(x,F)$.
Then $d(x,F)\leq(n+1)(4N+2)\lambda^{-1/\beta}$,
and we apply the induction estimate with $k=n$ to obtain
\[
R_\lambda\mathbf{1}_F(x)
\leq \frac{1}{2^n\lambda}
\leq \frac{2}{\lambda}
\exp\left[-\frac{\log 2}{4N+2}
d(x,F)\lambda^{1/\beta}\right].
\]

For general $F$,
apply the estimate to an increasing sequence of finite-measure Borel subsets
with union $F$.
Their distances from $x$ are at least $d(x,F)$,
so monotone convergence theorem proves the desired result.
\end{proof}

\begin{proposition}\label{prop:analytic-semigroup-decay}
There are $C,c>0$ such that, for  every Borel set $F\subset M$
and almost every $x\in M$,
\[
  P_t\mathbf{1}_F(x)
  \leq C\exp\left[
  -c\left(\frac{d(x,F)}{t^{1/\beta}}\right)^{\beta/(\beta-1)}
  \right],\qquad 0<t<\delta r_0^\beta.
\]
\end{proposition}

\begin{proof}
For $\lambda>0$ with $\lambda\geq\lambda_0/(\delta r_0^\beta)$,
put
\[
F_\lambda=\{y\in M:d(y,F)<\lambda^{-1/\beta}\}.
\]
If $z\in F$ and $0<s<1/\lambda$,
then $B(z,\delta s^{1/\beta})\subset F_\lambda$.
Integrating \eqref{eq:nle} on this ball and using \eqref{eq:vol},
we obtain $P_s\mathbf{1}_{F_\lambda}(z)\geq c_1c_2\delta^\alpha$
for almost every $z\in F$ at each such $s$.
It follows from \eqref{eq:analytic-resolvent-representation} that for almost every $z\in F$
\[
R_\lambda\mathbf{1}_{F_\lambda}(z)
\geq c_1c_2\delta^\alpha\int_0^{1/\lambda}e^{-\lambda s}\,\mathrm{d}s
\geq \frac{c_1c_2\delta^\alpha}{2\lambda}.
\]
Consequently,
$\mathbf{1}_F\leq (2\lambda/(c_1c_2\delta^\alpha))
R_\lambda\mathbf{1}_{F_\lambda}$ almost everywhere on $M$.

It holds by \eqref{eq:analytic-resolvent-representation} and the semigroup property that
\[
e^{-\lambda t}P_t R_\lambda\mathbf{1}_{F_\lambda}
=\int_t^\infty e^{-\lambda s}P_s\mathbf{1}_{F_\lambda}\,\mathrm{d}s
\leq R_\lambda\mathbf{1}_{F_\lambda}.
\]
We deduce from Lemma~\ref{lem:resolvent-decay} that for almost every $x\in M$
\begin{align*}
P_t\mathbf{1}_F(x)
\leq \frac{2\lambda}{c_1c_2\delta^\alpha}
e^{\lambda t}R_\lambda\mathbf{1}_{F_\lambda}(x)
\leq C\exp\left[\lambda t-c\,d(x,F_\lambda)\lambda^{1/\beta}\right]
\leq C\exp\left[\lambda t-c\,d(x,F)\lambda^{1/\beta}\right],
\end{align*}
where we used $d(x,F_\lambda)\geq d(x,F)-\lambda^{-1/\beta}$
and absorbed the resulting constant into $C$.

Choosing
\[
\lambda=\lambda_*=
\left(\frac{cd(x,F)}{\beta t}\right)^{
\frac{\beta}{\beta-1}}
\]
when $\lambda_*\geq\lambda_0/(\delta r_0^\beta)$,
and using the sub-Markov property to enlarge $C$
when $\lambda_*<\lambda_0/(\delta r_0^\beta)$,
we obtain the desired estimate.

The case $d(x,F)=0$ also follows directly from the sub-Markov property.
Thus the desired result holds for
almost every $x\in M$.
\end{proof}

\subsection{An analytic proof of the main result}

\begin{proof}[Analytic proof of Theorem~\ref{thm:main}]
Choose nonempty open subsets $U\ni x,V\ni y$ with $d(U,V)>0$ and put
\[
\mathcal{U}=\left\{z\in M:
d(z,U)<\frac{d(U,V)}{2}\right\}.
\]
Then $d(x,M\setminus \mathcal{U})\geq d(U,V)/2$ and $d(y,\mathcal{U})\geq d(U,V)/2$.

Splitting the semigroup integral over $M\setminus\mathcal{U}$ and $\mathcal{U}$,
and using symmetry and \eqref{eq:due} from \cite[Theorem~3.3]{ghl08},
we obtain for almost every $(x,y)\in U\times V$
\begin{align*}
p_t(x,y)
&=\int_M p_{t/2}(x,z)p_{t/2}(z,y)\,\mu(\mathrm{d}z)\\
&\leq Ct^{-\alpha/\beta}\left[
\int_{M\setminus \mathcal{U}}p_{t/2}(x,z)\,\mu(\mathrm{d}z)
+\int_\mathcal{U}p_{t/2}(y,z)\,\mu(\mathrm{d}z)
\right]\\
&\leq Ct^{-\alpha/\beta}\exp\left[-c\left(\frac{d(x,M\setminus \mathcal{U})^\beta}{t/2}\right)^{1/(\beta-1)}\right]
 +Ct^{-\alpha/\beta}\exp\left[-c\left(\frac{d(y,\mathcal{U} )^\beta}{t/2}\right)^{1/(\beta-1)}\right]\\
&\leq Ct^{-\alpha/\beta}
\exp\left[
-c\left(\frac{d(U,V)^\beta}{t}
\right)^{1/(\beta-1)}
\right]
\end{align*}
where the last inequality follows from Proposition~\ref{prop:analytic-semigroup-decay}.

Taking $U,V$ from a countable base of $M$ and
shrinking the corresponding neighborhoods to
$x$ and $y$, respectively, yields
\[
p_t(x,y)\leq Ct^{-\alpha/\beta}
\exp\left[
-c\left(\frac{d(x,y)}{t^{1/\beta}}\right)^{\beta/(\beta-1)}
\right]
\]
for almost every $x,y\in M$.
This completes the proof.
\end{proof}

\section{Sub-Gaussian Upper Bounds via Probabilistic Methods}\label{sec:upper-prob}

The probabilistic proof is based on an estimate for rapid radial passages.
The key step is to use path continuity and the near-diagonal lower estimate
to bound the probability of crossing a sequence of disjoint annuli in a short time.
We then iterate this estimate on deterministic time grids to obtain
an exponential bound for the probability of reaching a given radius by time $t$,
and combine it with the on-diagonal upper estimate to obtain
the full sub-Gaussian upper bound.
\subsection{Rapid radial passages}\label{sec:passages}
Let $(X_s)_{s\ge0}$ be the associated Hunt process with lifetime $\zeta$.
Locality gives continuous paths on $[0,\zeta)$ for almost every starting point \cite[Theorems~7.2.1 and 7.2.2]{fot94}.
Write $Q_s(x,A)=\mathbb P_x(X_s\in A)$ for the transition measure on $M$.
At each fixed $s>0$,
association with $P_s$ yields
\begin{equation}\label{eq:row}
  Q_s(x,\mathrm{d}y)=p_s(x,y)\,\mu(\mathrm{d}y)\quad\text{on }M
\end{equation}
for almost every $x\in M$.
Here equality holds as finite measures,
with a single exceptional set for all Borel target sets at that time.
To see this,
use association first for indicators from a countable determining class of relatively compact Borel sets.
All finite unions of elements of a countable base of relatively compact
open sets form such a class and include an increasing exhaustion of $M$.
These indicators belong to $L^2$.
Sub-Markovianity gives finiteness of the heat kernel rows by the exhaustion.
Equality then extends to the Borel $\sigma$-algebra by uniqueness of finite measures.
For $R>0$, let $\sigma_R=\inf\{0\le s<\zeta:d(X_s,X_0)\ge R\}$ with the convention $\inf\varnothing=\infty$.
For an open set $\Omega\subset M$,
define the exit time $\tau_\Omega=\inf\{s>0:X_s\notin\Omega\}$ with the cemetery state outside $\Omega$.
Under $\mathbb P_x$, we have $X_0=x$. Since the process is at the cemetery state from $\zeta$ onward,
we obtain $\tau_{B(x,R)}=\min(\sigma_R,\zeta)$ and hence
\begin{align*}
  \sigma_R=
  \begin{cases}
    \tau_{B(x,R)},
    & \tau_{B(x,R)}<\zeta,\\
    \infty,
    & \tau_{B(x,R)}=\zeta .
  \end{cases}
\end{align*}
On $\{\sigma_R<\infty\}$,  the passage occurs before the lifetime and $d(X_{\sigma_R},x)=R$ by continuity.

\begin{lemma}\label{lem:rapid}
There is $\varepsilon>0$ such that for every $r>0$ and almost every $x\in M$,
\[
  \mathbb P_x(\sigma_r\le t)\le 1/4,
  \qquad t=\varepsilon\min(r,r_0)^\beta .
\]
\end{lemma}

\begin{proof}
Integrating \eqref{eq:nle} and using \eqref{eq:vol},
we obtain for each fixed $0<s<\delta r_0^\beta$ and almost every $z$
\begin{equation}\label{eq:mass}
  Q_s(z,B(z,\delta s^{1/\beta}))
  =\int_{B(z,\delta s^{1/\beta})}
  p_s(z,y)\,\mu(\mathrm{d}y)
  \geq c_2s^{-\alpha/\beta}V(z,\delta s^{1/\beta})
  \geq c_1c_2\delta^\alpha>0,
\end{equation}
where the equality follows from \eqref{eq:row}.
First choose a sufficiently large integer $N$ and then a sufficiently small
$\varepsilon>0$, so that
\[
c_1c_2\delta^\alpha N\geq 4,
\qquad 2\varepsilon<\delta,
\qquad (4N+2)\varepsilon^{1/\beta}<1.
\]
Then $2t<\delta r_0^\beta$ and $(4N+2)t^{1/\beta}<r$.
For $x\in M$ and $1\leq j\leq N$,
define the open annuli centered at $x$ by
\begin{align*}
  A_j &= \{y\in M: (4j-2)t^{1/\beta}<d(x,y)<(4j+2)t^{1/\beta}\},\\
  B_j &= \{y\in M: (4j-1)t^{1/\beta}<d(x,y)<(4j+1)t^{1/\beta}\}.
\end{align*}
The annuli $A_j$ are pairwise disjoint.
If $z\in B_j$ and $0<s\leq t$,
the triangle inequality gives $B(z,\delta(2t-s)^{1/\beta})\subset A_j$,
because $\delta(2t-s)^{1/\beta}\le 2^{1/\beta}\delta t^{1/\beta}<t^{1/\beta}$
by $\delta\le1/2$ and $\beta>1$.

For each integer $m\geq1$,
let $G_m=(0,t]\cap t2^{-m}\mathbb N$ and let
$T_j^{(m)}$ be the first time $s\in G_m$ with $X_s\in B_j$,
with value $\infty$ if the grid misses $B_j$.
For $s\in G_m$, the event $\{T_j^{(m)}=s\}$ is measurable at time $s$,
namely,
\[
\{T_j^{(m)}=s\}
=\{X_s\in B_j\}
\cap\bigcap_{s'\in G_m:\ s'<s}\{X_{s'}\notin B_j\}
\in\mathcal{F}_s.
\]
Since $s\in G_m\subset(0,t]$, we have $t\le 2t-s<2t$.

By the Markov property, on $\{T_j^{(m)}=s\}$,
the conditional probability of $X_{2t}\in A_j$ given $\mathcal F_s$ is $Q_{2t-s}(X_s,A_j)$.
For almost every $x\in M$,
the law of $X_s$ restricted to $M$ is absolutely continuous by \eqref{eq:row},
so $X_s$ avoids the exceptional starting points for \eqref{eq:mass} at time $2t-s$.
The observation times and the time differences used here form countable sets.
Consequently, for almost every $x$ and all $m,j$,
\begin{align}\label{eq:annulus}
\mathbb P_x(X_{2t}\in A_j,\ T_j^{(m)}\le t)
&=\sum_{s\in G_m}\mathbb E_x\left[\mathbf 1_{\{T_j^{(m)}=s\}}\mathbf 1_{\{X_{2t}\in A_j\}}\right]\nonumber \\
&=\sum_{s\in G_m}\mathbb E_x\left[
    \mathbf 1_{\{T_j^{(m)}=s\}}\mathbb E_x[\mathbf 1_{\{X_{2t}\in A_j\}}\mid\mathcal F_s]
    \right]\nonumber \\
&=\sum_{s\in G_m}
    \int_{\{T_j^{(m)}=s\}}
    Q_{2t-s}(X_s,A_j)\,\mathrm d\mathbb P_x
    \nonumber\\
&\ge\sum_{s\in G_m}
    \int_{\{T_j^{(m)}=s\}}
    Q_{2t-s}\bigl(X_s,B(X_s,\delta(2t-s)^{1/\beta})\bigr)
    \,\mathrm d\mathbb P_x
    \nonumber\\
&\ge c_1c_2\delta^\alpha
    \sum_{s\in G_m}\mathbb P_x(T_j^{(m)}=s)
=c_1c_2\delta^\alpha
    \mathbb P_x(T_j^{(m)}\le t).
\end{align}
On $\{\sigma_r\le t\}$,
continuity forces the path to visit each open $B_j$
during a nonempty time interval before the passage.
The nested grids eventually sample each such interval.
Letting $m\to\infty$ in \eqref{eq:annulus},
we sum over $1\leq j\leq N$ and use the disjointness of the annuli $A_j$ to obtain
\[
\mathbb P_x(\sigma_r\le t)
\le\frac{1}{c_1c_2\delta^\alpha N}
\sum_{j=1}^N\mathbb P_x(X_{2t}\in A_j)
\le1/4,
\]
which proves the lemma.
\end{proof}

\begin{proposition}\label{prop:tail}
There are $C,c>0$ such that for each $R>0$ and almost every $x\in M$,
\[
  \mathbb P_x(\sigma_R\le t)\le C\exp\left[-c(R^\beta/t)^{1/(\beta-1)}\right],
  \qquad 0<t<r_0^\beta .
\]
\end{proposition}

\begin{proof}
Fix integers $n,m\geq1$ and set $r=R/(4n)$.
For $1\leq k\leq n$, define the recording times recursively on the grid
of step $2^{-m}$ by
\begin{align*}
  S_0^{(m)}&=0,\\
  S_{k}^{(m)}
  &=\inf\{s\in 2^{-m}\mathbb N:
    S_{k-1}^{(m)}<s<\zeta,\ d(X_s,X_{S_{k-1}^{(m)}})\ge r \}
    \quad\text{if }S_{k-1}^{(m)}<\infty,\\
  S_{k}^{(m)}&=\infty \quad\text{if }S_{k-1}^{(m)}=\infty .
\end{align*}
At each step,
the next recorded time is the first grid time before the lifetime
at which the process has moved at least $r$ from its last recorded position.
The empty infimum convention makes every subsequent recording time infinite once a step fails.

Set $\lambda=(\log 4)/t_0$, where $t_0=\varepsilon\min( r ,r_0)^\beta$ as in Lemma~\ref{lem:rapid}.
For almost every $z$,
the Borel function $g_m(z)=\mathbb E_z[e^{-\lambda S_1^{(m)}}]$
is at most $1/2$.
Indeed, since $S_1^{(m)}\ge\sigma_r$,
Lemma~\ref{lem:rapid} yields
\begin{align*}
\mathbb E_z\left[e^{-\lambda\sigma_ r }\right]
&=\mathbb E_z\left[e^{-\lambda\sigma_ r }
          \mathbf 1_{\{\sigma_ r \le t_0\}}\right]
 +\mathbb E_z\left[e^{-\lambda\sigma_ r }
          \mathbf 1_{\{\sigma_ r >t_0\}}\right]\\
&\le \mathbb P_z(\sigma_ r \le t_0)
    +e^{-\lambda t_0}\mathbb P_z(\sigma_ r >t_0)\le 1/4+e^{-\lambda t_0}
=1/2.
\end{align*}
For each $1\le k\le n$ and each finite value $s$ of $S_{k-1}^{(m)}$,
the Markov property at $s$ gives
\begin{align}\label{eq:iterate}
\mathbb{E}_x\left[e^{-\lambda S_{k}^{(m)}}\mathbf{1}_{\{S_{k-1}^{(m)}=s\}}\right]
&= \mathbb{E}_x\left[\mathbf{1}_{\{S_{k-1}^{(m)}=s\}}\mathbb{E}_x\left[e^{-\lambda S_{k}^{(m)}}\mid\mathcal{F}_s\right]\right]\nonumber\\
&= \mathbb E_x\left[\mathbf 1_{\{S_{k-1}^{(m)}=s\}}e^{-\lambda s}\mathbb E_{X_s}\left[e^{-\lambda S_1^{(m)}}\right] \right]\nonumber\\
&= \mathbb{E}_x\left[e^{-\lambda s}\mathbf{1}_{\{S_{k-1}^{(m)}=s\}}g_m(X_s)\right]
\le \frac{1}{2}\mathbb{E}_x\left[e^{-\lambda s}\mathbf{1}_{\{S_{k-1}^{(m)}=s\}}\right].
\end{align}
For $s>0$,
\eqref{eq:row} allows the bound on $g_m$ at $X_s$.
At $s=0$ the bound follows by choosing a starting point satisfying the estimate for $\sigma_ r $.
The Markov identity follows by considering finitely many future grid times first and then passing to the limit.
Summing \eqref{eq:iterate} over $s$ and iterating over $1\leq k\leq n$ yields
\[
\mathbb E_x \left[e^{-\lambda S_n^{(m)}}\right]\le \frac12\mathbb E_x \left[e^{-\lambda S_{n-1}^{(m)}}\right]\le \dots
\le \frac{1}{2^n}\mathbb E_x \left[e^{-\lambda S_0^{(m)}}\right]\le \frac{1}{2^n}\quad \text{for almost every } x \text{ and all }m.
\]

On $\{\sigma_R\le t\}$ put $s_m=2^{-m}\lfloor 2^m\sigma_R\rfloor$.
Consecutive grid positions up to $s_m$ are at distance less than $ r $ once the grid is sufficiently fine,
by uniform continuity on $[0,\sigma_R]$.
Each completed recorded displacement is then less than $2 r $,
and the displacement since the last recorded time is less than $ r $.

With fewer than $n$ recorded steps,
the triangle inequality would force $d(X_{s_m},x)\le(2n-1) r <R/2$.
The convergence $d(X_{s_m},x)\to R$ therefore implies $S_n^{(m)}\le s_m\le \sigma_R\le t$ for all sufficiently large $m$ on this event.
Fatou's lemma and the exponential moment bound now give
\begin{align}\label{eq:opt}
  \mathbb P_x(\sigma_R\le t)
  &\le\liminf_{m\to\infty}\mathbb P_x(S_n^{(m)}\le t)
   \le2^{-n}e^{\lambda t} \nonumber\\
  &\le \exp\left[K t(4n/R)^\beta+K t/r_0^\beta-n\log 2\right],
  \qquad K=(\log 4)/\varepsilon .
\end{align}
For each fixed $R$,
the exceptional sets may be taken simultaneously for all integers $n,m\geq1$.
The pathwise implication above then holds for every $t$ in the stated range.
Take $n=\lfloor \theta(R^\beta/t)^{1/(\beta-1)}\rfloor$ when $\theta(R^\beta/t)^{1/(\beta-1)}\ge2$,
where $\theta >0$ satisfies $K4^\beta \theta^{\beta-1}\le(\log 2)/2$.
The first term in the exponent of \eqref{eq:opt} is at most $n(\log 2)/2$,
and the term $Kt/r_0^\beta$ is bounded by $K$.
Since $n\ge (\theta/2)(R^\beta/t)^{1/(\beta-1)}$,
this proves the proposition for $\theta(R^\beta/t)^{1/(\beta-1)}\ge2$.
Increasing $C$ covers the remaining values by the bound $1$.
The proof is complete.
\end{proof}

\subsection{A probabilistic proof of the main result}
\begin{proof}[Probabilistic proof of Theorem~\ref{thm:main}]
For a fixed $R>0$,
\eqref{eq:row} yields
$$
\int_{d(x,z)\geq R}p_{t/2}(x,z)\,\mu(\mathrm{d}z)
=Q_{t/2}(x,M\setminus B(x,R))
=\mathbb P_x(X_{t/2}\in M\setminus B(x,R))
\leq \mathbb P_x(\sigma_R\leq t/2).$$
Apply Proposition~\ref{prop:tail} first at positive rational radii and take a common full-measure set of starting points.
For any positive real $R$,
the same bound for the integral follows by taking rational radii tending to $R$ from below.
This permits using $R=d(x,y)/2$ in the semigroup integral.

If $d(x,y)>0$,
each integration point $z$ is at distance at least $d(x,y)/2$ from $x$ or from $y$.
Splitting the semigroup integral over these two regions,
and using symmetry and \eqref{eq:due} from \cite[Theorem~3.3]{ghl08},
we obtain for almost every such $(x,y)\in M\times M$
\begin{align*}
  p_t(x,y)
  &=\int_M p_{t/2}(x,z)p_{t/2}(z,y)\,\mu(\mathrm{d}z)\\
  &\leq Ct^{-\alpha/\beta}\left[
  \int_{d(x,z)\ge d(x,y)/2}p_{t/2}(x,z)\,\mu(\mathrm{d}z)
  +\int_{d(y,z)\ge d(x,y)/2}p_{t/2}(y,z)\,\mu(\mathrm{d}z)
  \right]\\
  &\le C t^{-\alpha/\beta}
    \exp\left[-c\left(\frac{d(x,y)}{t^{1/\beta}}\right)^{\beta/(\beta-1)}\right],
\end{align*}
where the last inequality follows from Proposition~\ref{prop:tail}.
This completes the proof.
\end{proof}

{\footnotesize{\noindent \textbf{Acknowledgments.}
B. Li is supported by National Natural Science Foundation of China (12471094, 12571121), and Zhejiang NSF of China (LMS26A010015).
Y. Sun is supported by National Natural Science Foundation of China ( 12371206)
and the Fundamental Research Funds for the Central Universities, No. 050-63263078.
The authors acknowledge the use of ChatGPT in the preparation of this manuscript and assume full
responsibility for its content.
}}

\begingroup
\footnotesize
\setlength{\parskip}{0pt}
\setlength{\itemsep}{0pt}
\setlength{\parsep}{0pt}
\setlength{\topsep}{0pt}
\setlength{\partopsep}{0pt}
\setlength{\labelsep}{0.3em}
\renewcommand{\baselinestretch}{0.95}\selectfont
\emergencystretch=1em

\footnotesize{

\smallskip

\noindent Bo Li (bli@zjxu.edu.cn), Institute of Mathematics, Jiaxing University, Jiaxing 314001,  China.

\smallskip

}

\endgroup

\end{document}